# Boussinesq equation

Dana Seiova

Yerkezhan Assylbek

Akzhan Bakibayeva

Olzhas Akbayev

**Table of content**

# Introduction

In this project our team was assigned problem #7. The equation is as follows.

$$\frac{\partial^2 u}{\partial t^2} - \frac{\partial^2 u}{\partial x^2} - \frac{\partial^4 u}{\partial x^4} - 3\frac{\partial}{\partial x^2}u^2 = 0$$

In general, this equation is special case of Boussinesq equation, which is shown below (1). Therefore, in order reader to better undersand how it comes, we decided to concentrate on Boussinesq equation.

Boussinesq equation (1) belongs to Korteweg-de Vries kind of equations (Han & Yarkony, 2011). The equation describes the motion of long waves in two dimensions under the gravitation (Han & Yarkony, 2011). This equation has this form:

$$\frac{\partial^2 u}{\partial t^2} - c\frac{\partial^2 u}{\partial x^2} - \frac{\partial^4 u}{\partial x^4} - \frac{\partial^2}{\partial x^2}u^2 = 0. \qquad (1)$$

Here, we differentiate u = u(x, t) to the needed order. With c=-1 we have well-posed Boussinesq equation, while with c=1 not well-posed classical Boussinesq equation (Gradinaru, Hagedorn, & Joye, 2010). On the other hand, we have corrected equation:

$$\frac{\partial^2 u}{\partial t^2} - \frac{\partial^2 u}{\partial x^2} - \frac{\partial^4 u}{\partial x^4} - \frac{\partial^2}{\partial x^2}u^2 = 0 . \qquad (2)$$

In physical research of nonlinear waves spread in waveguides, it is necessary to take into account the interaction of waveguides with external environment (Berezhkoyskii, 2011). Therefore, the possibility of energy exchange through the outer surface of waveguides, is also should be taken into account (Berezhkoyskii, 2011). When the energy exchange between a rod and external environment occurs, it is important to estimate the deformation dispersion of wave in viscous environment (Han & Yarkony, 2011). The equation with two dispersion, which is derived from Hamilton's principle, has the third order:

$$\frac{\partial^2 u}{\partial t^2} - \frac{\partial^2 u}{\partial x^2} = \frac{1}{4}\left(cu^3 + 6u^2 + a\frac{\partial^2}{\partial t^2}u^2 - b\frac{\partial^2}{\partial x^2}u^2 + d\frac{\partial^2}{\partial t^2}u^2\right)\frac{\partial^2}{\partial x^2}$$

(3)

where a, b, c, d are any real positive numbers.

Symmetrical reductions and exact solutions have a lot of various important applications in the context of differential equations (Qin, 2006).

Since the solutions of partial differential equations are asymptotically approach to solutions of lower dimensional equations, which are obtained by symmetry reductions (Hsiao, Hu, & Hwung, 2010). Some of these special solutions describe important physical phenomena (Hsiao, Hu, & Hwung, 2010). In particular, exact solutions, derived by symmetrical methods, often effectively used in researches, such as asymptotic behavior and "blow-up" (Hsiao, Hu, & Hwung, 2010).

Moreover, explicit solutions (for instance, solutions, obtained by symmetrical methods) can play a significant role in construction and checking of numerical integrals (Yang, 2011).

In our research, we consider a generalized Boussinesq equation:

$$\Delta \equiv u_{tt} - u_{xxtt} + u_{xxxxtt} + cu_{xxxx} - u_{xx} = (f(u))_{xx} \qquad (4).$$

When $c \neq 0$, for some functions $f(u)$, the equation belongs to classical symmetry groups (Boulanba & Mellouk, 2014). By using symmetrical reduction, we reduce initial non-linear PDE to ordinary non-linear ODE. Exact solutions are then obtained from this equation by application of direct method (Shaoyong, 2014). Using the method of decomposition $\frac{G'}{G}$, which recently was introduced by Vang, we obtain new exact solutions to equation (4) for $f(u) = u^2$ (Shaoyong, 2014).

**Classical symmetry**

In order to apply the classical method Li to equation (4), let's consider infinitesimal transformations in a range of variables *x, t, u* of one-parameter Li group, given by the formula (Wazwaz, 2015) :

$$x^* = x + \epsilon \varepsilon(x, t, u) + O(\epsilon^2),$$

$$t^* = t + \epsilon \tau(x, t, u) + O(\epsilon^2), \quad (5)$$

$$u^* = u + \epsilon \eta(x, t, u) + O(\epsilon^2).$$

In these equations $\epsilon$ is a parameter of group transformation. We require that this transformation will leave invariant a set of equation solutions (4) (Quansen, Jiahong, & Wanrong, 2010). This requirement leads to over determined system of linear equations on small quantities *ξ(x, t, u), τ (x, t, u)* and *η(x, t, u)* (Quansen, Jiahong, & Wanrong, 2010*)*. Algebra associated with the infinitesimal transformation of symmetry is a set of vector fields of the form (Hu, Kukavica, & Ziane, 1997):

$$V = \varepsilon(x,t,u)\partial_x + \tau(x,t,u)\partial_t + \eta(x,t,u)\partial_u \qquad (6).$$

Functions $u = u(x,t)$ are invariant in relation to infinitesimal transformations v, essentially presents solutions to equation, obtained under "Invariant surface condition" (Zhang, & Lu, 2014).

$$\eta(x,t,u) - \varepsilon(x,t,u)u_x - \tau(x,t,u)u_t = 0 \qquad (7).$$

We can find symmetrical variables by solving this invariant surface condition (Hu, Kukavica, & Ziane, 1997). With reduction of PDE, we obtain ODE.

Let's move to classical analysis of symmetry group Lie for equation (4). The set of solutions for this equation is invariant in relation to transformation (5) (Yan, 2009). If $pr^{(6)} = v(\Delta) = 0$, where $\Delta=0$ and $pr^{(6)}v$ is continuity of the sixth order for the vector field (6) ) (Yan, 2009). Thus, there are appear sixteen variables that define equation for small quantities) (Yan, 2009). From this system, it follows that if $f$ is an arbitrary function, the only symmetry that equation (4) allows is symmetry (Yan, 2009):

$$\varepsilon = \lambda, \ \tau = \mu, \ \eta = 0, \qquad (8).$$

Given by a group of space-time transformations v1 = $\partial$x, v2 = $\partial$t (Himonas, & Mantzavinos, 2013). Substituting (8) in the invariant surface condition (7), we get a similar variable and such solution:

$$z = \mu x - \lambda t, u(x,t) = h(z) \qquad (9).$$

After substituting this ratio into (4), we have:

$$\lambda^2 \mu^4 h''''' + \mu^2(\mu^2 c - \lambda^2)h'''' + (\lambda^2 - \mu^2)h''' - \mu^2 f_h h'' - \mu^2 f''(h')^2 = 0 \qquad (10).$$

Integrating this equation two times in relation to $z$, we get:

$$\lambda^2 \mu^4 h'''' + \mu^2(\mu^2 c - \lambda^2)h'' + (\lambda^2 - \mu^2)h - \mu^2 f(h) + Az + B = 0 \qquad (11).$$

where A and B- are constants.

When, $f(u) = au + b$ the equation (4) allows additional symmetry to occur, and our equation became linear PDE, therefore we do not take into account this case (Himonas, & Mantzavinos, 2013).

# 3. Solving type of traveling wave

In this section we will get the solution in the form of traveling wave for the equation (11). We will use two procedures: direct method and method of G'/G expansion.

1. Direct Method

Consider equation (11) with A=B=0

$$h'''' + \frac{\mu^2 c - \lambda^2}{\mu^2 \lambda^2} h'' + \frac{\lambda^2 - \mu^2}{\mu^4 \lambda^2} h - \frac{1}{\mu^2 \lambda^2} f(h) = 0. \tag{12}$$

This equation can be rewritten in the form

$$h''''+bh''+F(h)=0 \tag{13}$$

where $b = (\mu^2 c - \lambda^2)/\mu^2 \lambda^2$

and

$$F(h) = \frac{\lambda^2 - \mu^2}{\mu^4 \lambda^2} h - \frac{1}{\mu^2 \lambda^2} f(h). \tag{14}$$

Equation (13) allows solutions of the form

$$h = \alpha H^\beta(z), \tag{15}$$

where α and β are parameters, and H(z) is solution of Jacobi's equation

$$(H')^2 = r + pH^2 + qH^4 \tag{16}$$

with constant r,p and q. Now we substitute (15) in (13) and by using identities for standard elliptic Jacobi's functions. We will get equation in terms of h and F(h). From these equations we will get function F(h), in which h is the solution of equation (13).

**Case 1**. If H(z) =sn z, then $h(z) = \alpha \, \text{sn}^\beta z$

Where cn is an elliptic function $u = \int_0^\phi \dfrac{d\theta}{\sqrt{1 - m\sin^2\theta}}$ and $\text{sn } u = \sin \phi$

$$F(h) = \alpha_1 h^{1+4/\beta} + \alpha_2 h^{1+2/\beta} + \alpha_3 h^{1-4/\beta} + \alpha_4 h^{1-2/\beta} + \alpha_5 h, \qquad (17)$$

Where

$$\begin{aligned}
\alpha_1 &= -\alpha^{-4/\beta}\left[\beta^4 m^4 + (6\beta^3 + 8\beta^2 + 4\beta)m^3 + (3\beta^2 + 2\beta)m^2\right], \\
\alpha_2 &= \alpha^{-2/\beta}\big[(2\beta^4 - 6\beta^3 + 8\beta^2 - 4\beta)m^4 + (12\beta^3 - 6\beta^2 + 8\beta)m^3 + \\
&\quad + (2\beta^4 + (6-b)\beta^2)m^2 + (6\beta^3 + 8\beta^2 + (4-b)\beta)m\big], \\
\alpha_3 &= -\alpha^{4/\beta}\beta^4 + 6\alpha^{4/\beta}\beta^3 - 11\alpha^{4/\beta}\beta^2 + 6\alpha^{4/\beta}\beta, \\
\alpha_4 &= (2\alpha^{2/\beta}\beta^4 - 12\alpha^{2/\beta}\beta^3 + 22\alpha^{2/\beta}\beta^2 - 12\alpha^{2/\beta}\beta)m^2 + \\
&\quad + (\alpha^{2/\beta}b - 4\alpha^{2/\beta})\beta + 2\alpha^{2/\beta}\beta^4 - 6\alpha^{2/\beta}\beta^3 + (8\alpha^{2/\beta} - \alpha^{2/\beta}b)\beta^2 + \\
&\quad + (6\alpha^{2/\beta}\beta^3 - 14\alpha^{2/\beta}\beta^2 + 8\alpha^{2/\beta}\beta)m, \\
\alpha_5 &= (-\beta^4 + 6\beta^3 - 11\beta^2 + 6\beta)m^4 + (-6\beta^3 + 14\beta^2 - 8\beta)m^3 + \\
&\quad + (-4\beta^4 + 12\beta^3 + (b-19)\beta^2 + (10-b)\beta)m^2 + \\
&\quad + (-12\beta^3 + 6\beta^2 + (b-8)\beta)m - \beta^4 + b\beta^2.
\end{aligned} \qquad (18)$$

**Case 2**. If H(z)=cn z, then $h(z) = \alpha \, \text{cn}^\beta z$

and

$$F(h) = \beta_1 h^{1+4/\beta} + \beta_2 h^{1+2/\beta} + \beta_3 h^{1-4/\beta} + \beta_4 h^{1-2/\beta} + \beta_5 h, \tag{19}$$

Where

$$\begin{aligned}
\beta_1 &= \alpha^{-4/\beta}\beta m^2(\beta^3 m^2 + 6\beta^2 m + 8\beta m + 4m + 3\beta + 2), \\
\beta_2 &= -\alpha^{-2/\beta}\beta m(4\beta^3 m^3 - 6\beta^2 m^3 + 8\beta m^3 - 4m^3 + 18\beta^2 m^2 + b + \\
&\quad + 2\beta m^2 + 12m^2 - 2\beta^3 m + b\beta m + 6\beta m - 6\beta^2 - 8\beta - 4), \\
\beta_3 &= \alpha^{4/\beta}(\beta - 3)(\beta - 2)(\beta - 1)\beta(m - 1)^2(m + 1)^2, \\
\beta_4 &= -\alpha^{2/\beta}(\beta - 1)\beta(m - 1)(m + 1)(4\beta^2 m^2 - 14\beta m^2 + 16m^2 + \\
&\quad + 6\beta m - 8m - 2\beta^2 + 4\beta + b - 4),
\end{aligned} \tag{20}$$

$$\begin{aligned}
\beta_5 &= \beta(6\beta^3 m^4 - 18\beta^2 m^4 + 27\beta m^4 - 14m^4 + 18\beta^2 m^3 - 20\beta m^3 + \\
&\quad + 16m^3 - 6\beta^3 m^2 + 12\beta^2 m^2 + 2b\beta m^2 - 13\beta m^2 - bm^2 + 6m^2 - \\
&\quad - 12\beta^2 m + 6\beta m + bm - 8m + \beta^3 - b\beta).
\end{aligned}$$

**Case 3**. If H(z) = dn z, then $h(z) = \alpha \operatorname{dn}^\beta z$

Where dn u is an elliptic function $\operatorname{dn} u = \sqrt{1 - m \sin^2 \phi}$.

and

$$F(h) = \gamma_1 h^{1-2/q} + \gamma_2 h + \gamma_3 h^{1+2/q} + \gamma_4 h^{1-4/q} + \gamma_5 h^{1+4/q}, \tag{21}$$

Where

$$\gamma_1 = -\alpha^{-4/\beta}m^{-4}\beta(8m^3 + 28\beta m^2 - 12m^2 + 12\beta^2 m - 28\beta m +$$
$$+ 16m + \beta^3 - 6\beta^2 + 11\beta - 6),$$
$$\gamma_2 = -\alpha^{-2/\beta}m^{-4}\beta(4m^5 + 28\beta m^4 - 12m^4 + 18\beta^2 m^3 - 42\beta m^3 -$$
$$- 2bm^3 + 16m^3 + 2\beta^3 m^2 - 12\beta^2 m^2 - b\beta m^2 - 34\beta m^2 + bm^2 +$$
$$+ 12m^2 - 36\beta^2 m + 84\beta m - 48m - 4\beta^3 + 24\beta^2 - 44\beta + 24),$$
$$\gamma_3 = -\alpha^{4/\beta}m^{-4}(\beta-3)(\beta-2)(\beta-1)\beta(m-1)^2(m+1)^2,$$

$$\gamma_4 = \alpha^{2/\beta}m^{-4}(\beta-1)\beta(m-1)(m+1)(6\beta m^3 - 8m^3 + 2\beta^2 m^2 - \quad (22)$$
$$- 10\beta m^2 - bm^2 + 12m^2 - 12\beta m + 16m - 4\beta^2 + 20\beta - 24),$$
$$\gamma_5 = -\beta(3\beta m^6 - 2m^6 + 6\beta^2 m^5 - 14\beta m^5 - bm^5 + 8m^5 + \beta^3 m^4 -$$
$$- 6\beta^2 m^4 - b\beta m^4 - 17\beta m^4 + bm^4 + 6m^4 - 36\beta^2 m^3 + 84\beta m^3 +$$
$$+ 2bm^3 - 48m^3 - 6\beta^3 m^2 + 36\beta^2 m^2 + 2b\beta m^2 - 38\beta m^2 -$$
$$- 2bm^2 + 24m^2 + 36\beta^2 m - 84\beta m + 48m + 6\beta^3 - 36\beta^2 +$$
$$+ 66\beta - 36)m^{-4}.$$

By substituting (17), (19) and (21) in (14), we will get 3 functions f(h), where h is the solution of (12), for cases 1,2 and 3 respectively.

**Case 1**. Let's substitute (17) in (14), as result we have

$$f(h) = \alpha'_1 h^{1+4/\beta} + \alpha'_2 h^{1+2/\beta} + \alpha'_3 h^{1-4/\beta} + \alpha'_4 h^{1-2/\beta} + \alpha'_5 h, \quad (23)$$

Where $\alpha'_i = -\mu^2\lambda^2\alpha_i$, i = 1, ..., 4, $\alpha'_5 = -\mu^2\lambda^2\alpha_5 + \mu^2(\lambda^2 - \mu^2)$,

and $\alpha_i$, i = 1,…5, are given by the formulas (18).

In this case the solution of equation (12) has form of $h(z) = \alpha\,\mathrm{sn}^\beta(z\mid m)$

and the exact solution of equation (4) is

$$u(x,t) = \alpha \operatorname{sn}^\beta(\mu x - \lambda t \mid m). \qquad (24)$$

where f(u) is got by substituting u instead of h in the formula (23).

Then we will provide some solutions, which are interesting from physical point of view.

- If m=0 and , µ = λ = $\sqrt{5/48}$ , α = 1 and β = 2 then we substitute these numbers in (23) then

$$f(h) = -\frac{5}{288}(12c^2 - 17)(2h - 1), \qquad (25)$$

and since sn(z | 0) = sin z, we have private solution

$$u(x,t) = \begin{cases} \sin^2\left[\sqrt{\dfrac{5}{48}}(x-t)\right], & |x-t| \leqslant \dfrac{2\pi}{k}, \\ 0, & |x-t| > \dfrac{2\pi}{k}, \end{cases} \qquad (26)$$

where k = $\sqrt{5/12}$ . This solution has two peaks. (Figure.1)

- If m = µ = β = 1, λ = 1/2 and α = ¼ then if we substitute if in (23)

$$f(h) = 1536h^5 + 32ch^3 - 168h^3 + \left(\frac{15}{4} - 2c\right)h, \qquad (27)$$

and since y sn(z | 1) = th z, and the solution of this equation (11) has form h(z) = (th z)/4.

In this case exact solution of equation (4), where f(u) is got by substituting u instead of h in the formula (27), thus

$$u(x,t) = \frac{1}{4}\operatorname{th}\left(x - \frac{t}{2}\right) \qquad (28)$$

and this represents the solution of kink-type.

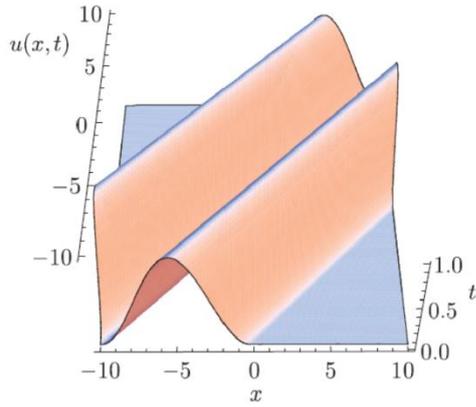
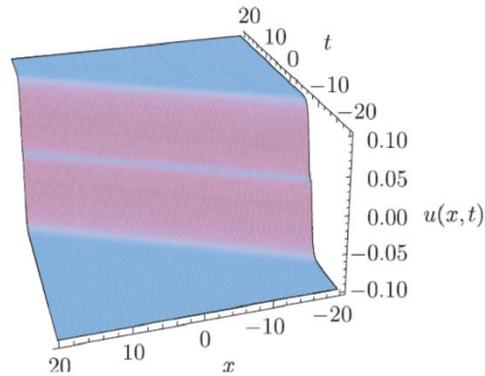

Figure.1                                      Figure.2

If we assume that $m = \mu = \alpha = 1$, $\lambda = 1/2$ and $\beta = 3$ equation (23) gives us

$$f(h) = 90h^{7/3} + 3(4c^2 - 69)h^{5/3} + \frac{3}{2}(4c^2 - 21)h^{1/3} - \frac{3}{4}(24c^2 - 197)h^{-1/3}, \qquad (29)$$

Then, the correct solution of equation (4), where f(u) is substituted instead of h in the formula (29), it turns to

$$u(x,t) = \text{th}^3\left(x - \frac{t}{2}\right) \qquad (30)$$

and presents by itself the solution of antikink (pic. 2)

**Case 2.**

By substituting (19) into (14), we get

$$f(h) = \beta'_1 h^{1+4/\beta} + \beta'_2 h^{1+2/\beta} + \beta'_3 h^{1-4/\beta} + \beta'_4 h^{1-2/\beta} + \beta'_5 h, \qquad (31)$$

Where

$$\beta'_i = -\mu^2 \lambda^2 \beta_i, \quad i = 1, \ldots, 4, \quad \beta'_5 = -\mu^2 \lambda^2 \beta_5 + \mu^2(\lambda^2 - \mu^2)$$

and $\beta_i$, for $i = 1, \ldots, 5$, are taken from the formula (20). The solution of equation (12) will be

$$h(z) = \alpha \, \mathrm{cn}^\beta(z \mid m).$$

Therefore, exact solution of equation (4), where we will get f(u) by substitution of u instead of h in equation (31), is

$$u(x,t) = \alpha \, \mathrm{cn}^\beta(\mu x - \lambda t \mid m). \tag{32}$$

By putting $m = 0$, $\mu = \lambda = \sqrt{\frac{5}{48}}$, $\alpha = 1$ and $\beta = 2$ in equation (31) we get

$$f(h) = \frac{5}{288}(12c^2 - 17)(2h - 1). \tag{33}$$

Because $\mathrm{cn}(z \mid 0) = \cos z$, the solution will be

$$u(x,t) = \begin{cases} \cos^2(\mu x - \lambda t), & |x - t| \leq \frac{\pi}{k}, \\ 0, & |x - t| > \frac{\pi}{k}, \end{cases} \tag{34}$$

where $k = \sqrt{\frac{5}{12}}$. It will be compact solution with one peak.

**Case 3.**

By substituting (21) into (14), we will get

$$f(h) = \gamma'_1 h^{1+4/\beta} + \gamma'_2 h^{1+2/\beta} + \gamma'_3 h^{1-4/\beta} + \gamma'_4 h^{1-2/\beta} + \gamma'_5 h, \tag{35}$$

where $\gamma'_i = -\mu^2 \lambda^2 \gamma_i, \; i = 1, \ldots, 4, \; \gamma'_5 = -\mu^2 \lambda^2 \gamma_5 + \mu^2(\lambda^2 - \mu^2)$, and $\gamma_i, \; i = 1, \ldots, 5,$ comes from (22). the solution of equation (12) has $h(z) = \alpha \, \mathrm{dn}^\beta(z \mid m)$. Therefore the exact solution of equation (4), where we will get f(u) by substitution of u instead of h in equation (35), is

$$u(x,t) = \alpha \, \text{dn}^\beta(\mu x - \lambda t \mid m). \tag{36}$$

By putting m = λ = μ = α = 1 and β = 2 in equation (35), we will get

$$f(h) = 120h^3 - 6(c^2 + 19)h^2 + 4(c^2 + 3)h. \tag{37}$$

Because

$$\text{dn}(z \mid 1) = \text{sech}\, z$$

we can the partial solution of $h(z) = \text{sech}^2 z$. Therefore, the exact solution of equation (4), where f(u) we get from substitution $u$ instead of $h$ to the equation (37), and we'll get

$$u(x,t) = \text{sech}^2(x - t) \tag{38}$$

and presents soliton, moving along straight line with the constant velocity.

**G′ /G-expansion method.**

Let's look at generalized Boussinesq equation (4) where $f(u) = u^2$

and take his solutions for traveling waves. In this case his equation will take the form of

$$\mu^4\lambda^2 h'''' + \mu^2(\mu^2 c - \lambda^2)h'' + (\lambda^2 - \mu^2)h - \mu^2 h^2 + Az + B = 0. \tag{39}$$

to apply G′ /G-expansion method to this equation, assume that A=0 and the solutions could be presented as polynomial of G′ /G:

$$h = \sum_{i=0}^{n} a_i \left(\frac{G'}{G}\right)^i, \qquad (40)$$

where $a_i$ for, $i = 0, \ldots, n$, - some constants, and $a_i \neq 0$, and $G = G(z)$, and satisfy second order linear differential equations

$$G''(z) + \alpha G(z) + \beta G = 0, \qquad (41)$$

$\alpha$ and $\beta$ will be defined further.

General solutions of equation (41) comes to the following:

If $\alpha^2 - 4\beta > 0$, then

$$G(z) = c_1 \operatorname{ch}\left(\frac{z\alpha}{2} - \frac{1}{2} z \sqrt{\alpha^2 - 4\beta}\right) + c_2 \operatorname{ch}\left(\frac{\alpha z}{2} + \frac{1}{2} \sqrt{\alpha^2 - 4\beta} z\right) -$$
$$- c_1 \operatorname{sh}\left(\frac{z\alpha}{2} - \frac{1}{2} z \sqrt{\alpha^2 - 4\beta}\right) - c_2 \operatorname{sh}\left(\frac{\alpha z}{2} + \frac{1}{2} \sqrt{\alpha^2 - 4\beta} z\right); \qquad (42)$$

If $\alpha^2 - 4\beta < 0$, then

$$G(z) = \left(c_2 \cos\left(\frac{1}{2} z \sqrt{4\beta - \alpha^2}\right) + c_1 \sin\left(\frac{1}{2} z \sqrt{4\beta - \alpha^2}\right)\right)\left(\operatorname{ch}\left(\frac{z\alpha}{2}\right) - \operatorname{sh}\left(\frac{z\alpha}{2}\right)\right); \qquad (43)$$

If $\alpha^2 = 4\beta$, then

$$G(z) = (c_2 + c_1 z)\left(\operatorname{ch}\left(\frac{z\alpha}{2}\right) - \operatorname{sh}\left(\frac{z\alpha}{2}\right)\right). \qquad (44)$$

By using the relation of (40) and (41), we will get

$$h^2 = a_n^2 \left(\frac{G'}{G}\right)^{2n} + \cdots,$$
$$h'''' = n(n+1)(n+2)(n+3)a_n \left(\frac{G'}{G}\right)^{n+4} + \cdots. \qquad (45)$$

Considering summand $h''''$ and $h^2$ in equation (39), from the formulas (45) we require $n+4=2n$, means $n=4$, then we can rewrite equation (40) as

$$h = a_0 + a_1\left(\frac{G'}{G}\right) + a_2\left(\frac{G'}{G}\right)^2 + a_3\left(\frac{G'}{G}\right)^3 + a_4\left(\frac{G'}{G}\right)^4, \qquad (46)$$

and $a_4 \neq 0$.

By substituting general solution of equation (41) into (46), we will get the following: for the solution, by formula (42)

$$\begin{aligned}h_1(z) &= \frac{a_4\alpha^4}{16} - \frac{a_3\alpha^3}{8} + \frac{a_2\alpha^2}{4} - \frac{a_1\alpha}{2} + a_0 + \\ &+ \frac{\sqrt{\alpha^2 - 4\beta}}{8}(-2a_4\alpha^3 + 3a_3\alpha^2 - 4a_2\alpha + 4a_1)H_1 + \\ &+ \frac{\alpha^2 - 4\beta}{8}(3a_4\alpha^2 - 3a_3\alpha + 2a_2)(H_1)^2 + \\ &+ \frac{(\alpha^2 - 4\beta)^{3/2}}{8}(a_3 - 2a_4\alpha)(H_1)^3 + \frac{(\alpha^2 - 4\beta)^2}{16}a_4(H_1)^4,\end{aligned} \qquad (47)$$

where

$$H_1(z) = \frac{c_2\,\text{ch}(\sqrt{\alpha^2 - 4\beta}\,z/2) + c_1\,\text{sh}(\sqrt{\alpha^2 - 4\beta}\,z/2)}{c_1\,\text{ch}(\sqrt{\alpha^2 - 4\beta}\,z/2) + c_2\,\text{sh}(\sqrt{\alpha^2 - 4\beta}\,z/2)};$$

For the solution, by formula (43),

$$\begin{aligned}h_2(z) &= \frac{a_4\alpha^4}{16} - \frac{a_3\alpha^3}{8} + \frac{a_2\alpha^2}{4} - \frac{a_1\alpha}{2} + a_0 + \frac{a_4}{16}(4\beta - \alpha^2)^2(H_2)^4 + \\ &+ \frac{a_3 + 2a_4\alpha}{8}(4\beta - \alpha^2)^{3/2}H_2 + \frac{3a_4\alpha^2 - 3a_3\alpha + 2a_2}{8}(4\beta - \alpha^2)(H_2)^2 + \\ &+ \frac{-2a_4\alpha^3 + 3a_3\alpha^2 - 4a_2\alpha + 4a_1}{8}\sqrt{4\beta - \alpha^2}\,(H_2)^3,\end{aligned} \qquad (48)$$

where

$$H_2(z) = \frac{c_2\cos(\sqrt{4\beta - \alpha^2}\,z/2) + c_1\sin(\sqrt{4\beta - \alpha^2}\,z/2)}{c_1\cos(\sqrt{4\beta - \alpha^2}\,z/2) - c_2\sin(\sqrt{4\beta - \alpha^2}\,z/2)};$$

For the solution, by formula (44),

$$h_3(z) = \frac{a_4\alpha^4}{16} - \frac{a_3\alpha^3}{8} + \frac{a_2\alpha^2}{4} - \frac{a_1\alpha}{2} + a_0 + \frac{a_4 c_1^4}{(c_1 z + c_2)^4} + \frac{(a_3 - 2a_4\alpha)c_1^3}{(c_1 z + c_2)^3} +$$
$$+ \frac{(3a_4\alpha^2 - 3a_3\alpha + 2a_2)c_1^2}{2(c_1 z + c_2)^2} + \frac{(-2a_4\alpha^3 + 3a_3\alpha^2 - 4a_2\alpha + 4a_1)c_1}{4(c_1 z + c_2)}. \quad (49)$$

Then we identify $a_i$, $i = 0, \ldots, 4$. From formula (46) we will find $h^2 h''$ and $h''''$ and substitute expressions of these values to the formula (39). By equating coefficients of degrees of $\left(\frac{G'}{G}\right)^i$, $i = 0, \ldots, 4$, we will get system of algebraic equations, for $a_i$, α, β, λ, μ and B. the solution of which should be identified by the following formulas:

$$a_0 = \frac{1}{338\lambda^2\mu^2} \left[ 3(5915\alpha^4\lambda^4 + 910c\alpha^2\lambda^2 + 23c^2)\mu^4 - \right.$$
$$\left. - (169 + 6(455\alpha^2\lambda^2 + 23c))\mu^2\lambda^2 + 238\lambda^4 \right],$$

$$a_1 = \frac{420}{13}\alpha\left[(13\alpha^2\lambda^2 + c)\mu^2 - \lambda^2\right], \quad a_2 = \frac{420}{13}\left[(39\alpha^2\lambda^2 + c)\mu^2 - \lambda^2\right],$$
$$a_3 = 1680\alpha\lambda^2\mu^2, \quad a_4 = 840\lambda^2\mu^2, \quad (50)$$
$$\beta = \frac{1}{52\lambda^2\mu^2}\left[13\alpha^2\mu^2\lambda^2 - \lambda^2 + c\mu^2\right],$$
$$B = \frac{1}{114244\lambda^4\mu^2}\left[-133\lambda^4 - (72c - 169)\mu^2\lambda^2 + 36c^2\mu^4\right]\left[205\lambda^4 - \right.$$
$$\left. - (72c + 169)\mu^2\lambda^2 + 36c^2\mu^4\right].$$

By replacing this values into (47)-(49) we will get three types of solutions of traveling waves for generalized Boussinesq equation (4) with $f(u) = u^2$

Firstly,

$$u_1(x,t) = \frac{215040}{169\lambda^2\mu^2}(c\mu^2 - \lambda^2)^2 F_1^4 - \frac{6720}{169\lambda^2\mu^2}(c\mu^2 - \lambda^2)^2 F_1^2 +$$

$$+ \frac{69c^2\mu^2}{338\lambda^2} - \frac{69c}{169} + \frac{119\lambda^2}{169\mu^2} - \frac{1}{2}, \qquad (51)$$

where

$$F_1 = \frac{c_2 \operatorname{ch}[\omega_1(x\mu - t\lambda)] + c_1 \operatorname{sh}[\omega_1(x\mu - t\lambda)]}{c_1 \operatorname{ch}[\omega_1(x\mu - t\lambda)] + c_2 \operatorname{sh}[\omega_1(x\mu - t\lambda)]},$$

$$\omega_1 = \frac{1}{2}\sqrt{\frac{1}{13}\left(-\frac{c}{\lambda^2} + \frac{1}{\mu^2}\right)}.$$

Secondly, the solution of type

$$u_2(x,t) = \frac{105}{338\lambda^2\mu^2}(c\mu^2 - \lambda^2)^2 F_2^4 + \frac{105}{169\lambda^2\mu^2}(c\mu^2 - \lambda^2)^2 F_2^2 +$$

$$+ \frac{69c^2\mu^2}{338\lambda^2} - \frac{69c}{169} + \frac{119\lambda^2}{169\mu^2} - \frac{1}{2}, \qquad (52)$$

Where

$$F_2 = -\frac{c_2 \cos(\omega_2(x\mu - t\lambda)) + c_1 \sin(\omega_2(x\mu - t\lambda))}{c_2 \sin(\omega_2(x\mu - t\lambda)) - c_1 \cos(\omega_2(x\mu - t\lambda))},$$

$$\omega_2 = \frac{1}{2\sqrt{13}}\sqrt{\frac{c}{\lambda^2} - \frac{1}{\mu^2}}.$$

Finally, the solution is

$$u_3(x,t) = \frac{840\lambda^2\mu^2 c_1^4}{(c_2 - c_1 t\lambda + c_1 x\mu)^4} - \frac{420(\lambda^2 - c\mu^2)c_1^2}{13(-c_2 + c_1 t\lambda - c_1 x\mu)^2} + \frac{69c^2\mu^2}{338\lambda^2} - \frac{69c}{169} + \frac{119\lambda^2}{169\mu^2} - \frac{1}{2}. \qquad (53)$$

Let's mark that we got general expression for the equation (46), which was not derived in the work (21). By substituting $a_3 = a_4 = 0$, $a_1 = -12\delta\alpha$, $a_2 = -12\delta$ and $z = x - (a_0 + 8\delta\beta + \delta\alpha^2)t$ to the formula (47)-(49), we will get the solution of Korteweg-de Vries equation, which were mentioned in the work of Shaoyong (2014).

Given equation:

$d^2u/dt^2 - d^2u/dx^2 - d^4u/dx^4 - 3 d(u^2)/dx^2 = 0$

this can be written as:

$u_{tt} - u_{xx} - u_{xxxx} - 3(u^2)_{xx} = 0$ (1)

what we know is that this equation can be solved by direct method or by $G'/G$ expansion method.

Let's introduce new variable $z$;

For some $\mu$ and real valued $\lambda$:

$Z = \mu x - \lambda t$; $u(x,t) = h(z)$

So this can be done, because u is a function of x and t. then our equation (1) will be rewritten as:

$(\lambda^2 - \mu^2) \ast h'' + \mu^4 \ast h'''' - 3\mu^2 \ast ((h')^2)'' = 0$ (2)

Mention, that this equation is in terms of $z$. then by simple double integration by $z$, we get the following:

$(\lambda^2 - \mu^2) \ast h + \mu^4 \ast h'' - 3\mu^2 \ast h^2 + Az + B = 0$ (3)

Where A,B are constants. By putting A=B=0. Then,

$h'' + (\lambda^2 - \mu^2)/\mu^4 \ast h - 3\mu^2 \ast h^2/\mu^4 = 0$

$F(h) = (\lambda^2 - \mu^2)/\mu^4 \ast h - 3h^2/\mu^2$ (3a)

This equation allow solutions of the form $h = \alpha H^\beta(z)$, where α, β are parameters and $H(z)$- Jacobian solution.

$(H')^2 = r + pH^2 + qH^4$ with $p, r, q$- constants.

We have 3 cases:

Case 1:

If $H(z) = sn\ z$; then $h(z) = \alpha sn^\beta(z)$

And

$F(h) = \alpha_1 h^{1+4/\beta} + \alpha_2 h^{1+2/\beta} + \alpha_3 h^{1-4/\beta} + \alpha_4 h^{1-2/\beta} + \alpha_5 h$ (4)

Case 2:

If $H(z) = cn\ z$; then $h(z) = \alpha cn^\beta(z)$

And

$$F(h) = \beta_1 h^{1+4/\beta} + \beta_2 h^{1+2/\beta} + \beta_3 h^{1-4/\beta} + \beta_4 h^{1-2/\beta} + \beta_5 h \qquad (5)$$

Case 3:

If **H(z) = dn z**; then **h(z)= αdn$^\beta$(z)**

And

$$F(h) = \gamma_1 h^{1-4/q} + \gamma_2 h + \gamma_3 h^{1+2/q} + \gamma_4 h^{1-4/q} + \gamma_5 h^{1+4/q} \qquad (6)$$

After direct putting (4), (5), (6) in (3a) we get:

Case 1:

(4) in (3a):

$$f(h) = \alpha'_1 h^{1+4/\beta} + \alpha'_2 h^{1+2/\beta} + \alpha'_3 h^{1-4/\beta} + \alpha'_4 h^{1-2/\beta} + \alpha'_5 h$$

so **u(x,t) = αsn$^\beta$(μx- λ+m)**

Case 2:

(5) in (3a):

$$f(h) = \beta'_1 h^{1+4/\beta} + \beta'_2 h^{1+2/\beta} + \beta'_3 h^{1-4/\beta} + \beta'_4 h^{1-2/\beta} + \beta'_5 h$$

so **u(x,t) = αcn$^\beta$(μx- λ+m)**

Case 3:

(6) in (3a):

$$f(h) = \gamma'_1 h^{1+4/\beta} + \gamma'_2 h^{1+2/\beta} + \gamma'_3 h^{1-4/\beta} + \gamma'_4 h^{1-2/\beta} + \gamma'_5 h$$

so **u(x,t) = αdn$^\beta$(μx- λ+m)**

**Conclusion**

In this paper was shown solution for partial differential equation

$$\frac{\partial^2 u}{\partial t^2} - \frac{\partial^2 u}{\partial x^2} - \frac{\partial^4 u}{\partial x^4} - 3\frac{\partial}{\partial x^2} u^2 = 0$$

While trying to find solution for this equation, we studied Boussinesq equation (4), which allows two-parameter group symmetry, which satisfies for two-parameter group symmetries. To get the exact solution can be used two methods: direct and G′ /G-expansion method. With the help of direct method it is possible to figure out soliton type, kink and antikink type exact solutions. By G′ /G-expansion method where $f(u) = u^2$ is used for finding three new types of solutions of traveling waves.